\magnification=\magstep1   
\input amstex
\UseAMSsymbols
\input pictex 
\vsize=23truecm
\NoBlackBoxes
\parindent=18pt
  
   \font\rmk=cmr8


\def\op{\text{\rm op}}

\def\mod{\operatorname{mod}}

\def\Hom{\operatorname{Hom}}

\def\End{\operatorname{End}}
\def\Ext{\operatorname{Ext}}

\def\rad{\operatorname{rad}}
\def\add{\operatorname{add}}
\def\Ker{\operatorname{Ker}}
\def\Cok{\operatorname{Cok}}

\def\Tr{\operatorname{Tr}}

\def\bdim{\operatorname{\bold{dim}}}

\def\top{\operatorname{top}}

  \def\ss{\ssize }
\def\arr#1#2{\arrow <1.5mm> [0.25,0.75] from #1 to #2}

\def\s{\hfill \square} 
 
\vglue3.5truecm
\centerline{\bf On modules $M$ such that}
	\smallskip
\centerline{\bf both $M$ and $M^*$ are semi-Gorenstein-projective.} 
                     \bigskip\medskip
\centerline{Claus Michael Ringel, Pu Zhang}
                \bigskip\medskip

\noindent {\narrower Abstract:  \rmk Let $\ss A$ be an artin algebra.
An $\ss A$-module $\ss M$ is semi-Gorenstein-projective provided that
$\ss\Ext^i(M,A) = 0$ for all $\ss i\ge 1.$ If $\ss M$ is Gorenstein-projective, then
both $\ss M$ and its $\ss A$-dual $\ss M^*$ are semi-Gorenstein projective. 
As we have shown recently, the converse is not true, 
thus answering a question raised by Avramov and Martsinkovsky.
The aim of the present note is to analyse in detail the modules $\ss M$ 
such that both $\ss M$ and $\ss M^*$ are semi-Gorenstein-projective.
	\medskip
\noindent 
\rmk Key words. Gorenstein-projective module, semi-Gorenstein-projective module, 
finitistic dimension conjecture, Nunke condition. 
	\medskip
\noindent 
Supported by NSFC 11971304.
\par}
	\bigskip
{\bf 1. Introduction.} 
	\medskip
Let $A$ be an artin algebra. The modules to be considered are usually 
left $A$-modules of finite length. Given a module $M$, let $M^* = \Hom(M,A)$ 
be its $A$-dual, and $\phi_M\:M \to M^{**}$ the canonical map from $M$ to $M^{**}.$
We will have to deal with complexes $P_\bullet = (P_i,f_i\:P_i\to P_{i-1})$ of
projective modules. Such a complex is said to be {\it minimal} provided the image of
$f_i$ is contained in the radical of $P_{i-1}$, for all $i\in \Bbb Z.$

A module is said to be
{\it reduced} provided it has no non-zero projective direct summand. 
The Main Theorem 2.1 asserts that the (isomorphism classes of the) reduced modules $M$ 
such that both $M$ and $M^*$ are semi-Gorenstein-projective 
correspond bijectively to the (isomorphism classes of) minimal 
complexes $P_\bullet$ of projective modules with $H_i(P_\bullet) = 0$
for $i\neq 0,-1,$ and such that the $A$-dual complex $P_\bullet^*$ is acyclic.
The essential idea is Lemma 2.2 which shows in which way 
the canonical map $\phi_M$ is related to a four term exact sequence
of projective right modules with $M^*$ being the image of the middle map.
Let us mention that $\Tr$ provides a bijection between the reduced
semi-Gorenstein-projective modules $M$ with also $M^*$ semi-Gorenstein-projective 
on the one hand, and
the reduced $\infty$-torsionfree right modules $Z$ 
with $\Omega^2 Z$ being semi-Gorenstein-projective, on the other hand, see 2.4.
These results are summarized by exhibiting
the complexes $P_\bullet$ and $P_\bullet^*$ with the modules $M, M^*, M^{**},$
as well as $\Tr M$ and $\Tr M^*$ being inserted, see 2.5.
	
Section 3 is devoted to two special situations. First, in 3.1, we
consider the case that $M$ is a semi-Gorenstein-projective module $M$ with also
$M^*$ semi-Gorenstein-projective such that $\phi_M$ is either an epimorphism or a 
monomorphism. In 3.2, we deal with semi-Gorenstein-projective modules $M$
such that $M^*$ is projective or even zero. 	
Whereas there do exist modules $M$ such that
both $M$ and $M^*$ are semi-Gorenstein-projective, but not Gorenstein-projective (see [RZ1]
and [RZ2]), it is not known whether we may have in addition that $\phi_M$ is either an
epimorphism or a monomorphism. 
Also, it is not known whether there exists a semi-Gorenstein-projective
module which is not projective, such that $M^*$ is projective. 

A module $M$ will be said to be a {\it Nunke} module provided that $M$ is 
semi-Gorenstein-projective and $M^* = 0.$ Note that an indecomposable semi-Gorenstein-projective module $M$ such that $\phi_M$ is an epimorphism and $M^*$ is projective, 
is either itself projective or else a Nunke module, 
see Proposition 3.4. The remaining parts of section 3 
are devoted to Nunke modules. 
There is the old conjecture (one of the classical homological conjectures) that 
the only Nunke module
is the zero module. This conjecture and similar ones are discussed in 3.5.
In 3.6 to 3.9 we try to analyze the special case of a simple injective module $S$ which is
semi-Gorenstein-projective (thus $S$ is either projective or a Nunke module).

In the final section 4, we consider local algebras, and, in particular, those 
with radical cube zero 
(the {\it short} local algebras). 
In [RZ1] and [RZ2], we have exhibited short local algebras with modules
$M$ such that both $M$ and $M^*$ are semi-Gorenstein-projective, whereas $M$
is not Gorenstein-projective. Here we show: Let $A$ be a local algebra
and $M$ a module such that both $M$ and $M^*$ are semi-Gorenstein-projective.
If $M^*$ is projective, then also $M$ is projective. If $A$ is short, and
$\phi_M$ is a monomorphism or an epimorphism, then $\phi_M$ is an isomorphism, thus
$M$ is Gorenstein-projective. In this way, we see that for short local algebras, 
there are no non-trivial
examples of modules which satisfy the conditions discussed in section 3.
	\bigskip\bigskip
{\bf 2. The Main Theorem.} 
	\medskip
{\bf 2.1. Main Theorem.} {\it The isomorphism classes of the reduced
modules $M$ such that both
$M$ and $M^*$ are semi-Gorenstein-projective correspond bijectively to the
isomorphism classes of minimal 
complexes $P_\bullet$ of projective modules with $H_i(P_\bullet) = 0$
for $i\neq 0,-1,$ and such that the $A$-dual complex $P_\bullet^*$ is acyclic, as follows: 

First, let $M$ be a reduced module such that
both $M, M^*$ are semi-Gorenstein-projective. Take minimal projective resolutions 
$$
 \cdots @>>> P_2 @>f_2>> P_1 @>f_1>> P_0 @>e>> M @>>> 0 \qquad \text{and} \qquad
 0 @<<< M^* @<q<< Q_0 @<d_1<< Q_1 @<d_2<< Q_2 @<<<\cdots.
$$
For $i < 0$, let $P_i = Q_{-i-1}^*$ and $f_i = d_{-i}^*\:P_{-i} \to P_{-i-1}$,
Finally, let $f_0 = q^*\phi_Me\:P_0 \to P_{-1}$. In this way, we obtain a minimal complex
$P_\bullet$ of projective modules
$$
\cdots @>>>
 P_1 @>{\,f_1\ }>> P_0 @>{\ f_0\ }>> P_{-1} @>f_{-1}>> P_{-2}  @>>> \cdots.
$$
with $H_i(P_\bullet) = 0$
for $i\neq 0,-1,$ and such that the $A$-dual complex $P_\bullet^*$ is acyclic. 
By construction, $M = \Cok f_1$ and $M^* = \Cok d_1 = \Cok f_{-1}^*$. Also,
$$
 H_{0}(P_\bullet) = \Ker \phi_M \quad\text{\it and}\quad 
 H_{-1}(P_\bullet) = \Cok \phi_M.
$$
	\medskip
Conversely, let \ $P_\bullet = (P_i,f_i\:P_i\to P_{i-1})_i$ 
be a minimal complex of projective modules 
with $H_i(P_\bullet) = 0$
for $i\neq 0,-1$, and such that 
the $A$-dual complex $P_\bullet^*$ is acyclic.
Let $M = \Cok f_1$. Then $M$ is reduced and
both $M$ and $M^*$ are semi-Gorenstein-projective. 
Actually, $M^* = \Cok f_{-1}^*$, $M^{**} = \Ker f_{-1}$.}
	\medskip 
The essential part of the proof is the following general lemma which shows in which way
exact sequences of projective modules are related to the canonical maps $\phi_M\:M \to M^{**}.$
	\medskip
{\bf 2.2. Lemma.} {\it Let
$$
 Q_{-2} @<d_{-1}<< Q_{-1} @<\ d_0\ << Q_{0} @<\ d_1\ << Q_{1} 
$$
be a sequence of projective right modules with composition zero. Let $e\:Q_1^* \to M$
be a cokernel of $d_{-1}^*$ and $c\:Q_0 \to N$ a cokernel of $d_1.$
The following conditions are equivalent:
\item{\rm (i)} The sequence is exact.
\item{\rm(ii)} There exists an isomorphism $\zeta\:N \to M^*$ with $d_0^* = c^*\zeta^*\phi_Me.$
\par}
	\medskip 
If the condition (ii) is satisfied, then $q = \zeta c\:
Q_0 \to M^*$ is a cokernel of $d_1$ and we have the following commutative diagrams:
$$
{\beginpicture
    \setcoordinatesystem units <1.5cm,1cm>
\put{\beginpicture
\put{$Q_{-1}$} at 0 0
\put{$M^*$} at 1 -1
\put{$Q_{0}$} at 2 0
\arr{0.8 -.8}{0.2 -0.2}
\arr{1.7 0}{0.3 0}
\arr{1.8 -.2}{1.2 -.8}
\put{$e^*$} at 0.4 -0.7
\put{$q$} at 1.55 -0.7
\put{$d_0$\strut} at 1 0.3

\endpicture} at 0 0.15

\put{\beginpicture
\put{$Q_{-1}^*$} at 0 0
\put{$M$} at 1 -1
\put{$M^{**}$} at 2 -1
\put{$Q_{0}^*$} at 3 0
\put{$e$} at 0.5 -0.7
\put{$\phi_M$} at 1.5 -1.23
\put{$q^*$} at 2.6 -0.7
\put{$d_0^*$} at 1.5 0.3
\arr{0.3 -0.3}{0.8 -.8}
\arr{0.3 0}{2.7 0}
\arr{1.2 -1}{1.7 -1}
\arr{1.2 -1.02}{1.7 -1.02}
\arr{1.2 -.98}{1.7 -.98}
\arr{2.2 -.8}{2.8 -.2}
\endpicture} at 3.5 0
\endpicture}
$$
	\medskip
Before we start with the proof, let us recall:
The exact sequence $0 @<<< N @<c<< Q_0 @<d_1<< Q_0$
yields the exact sequence $0 @>>> N^* @>c^*>> Q_0^* @>d_1^*>> Q_1^*$, thus
$\Ker d_1^* = (\Cok d_1)^*.$ 
Similarly, the exact sequence $Q_{-2}^* @>d_{-1}^*>> Q_{-1}^* @>e>> M @>>> 0$
yields the exact sequence $Q_{-2} @<d_{-1}<< Q_{-1} @<e^*<< M^* @<<< 0$, thus
$\Ker d_{-1} = (\Cok d_{-1}^*)^*.$
	\medskip
Proof of Lemma. (i) implies (ii). Since the sequence $Q_\bullet$ is exact, 
the cokernel $N$ of $d_1$ is equal to the kernel $M^*$ 
of $d_{-1}$, thus there is an isomorphism $\zeta\: N \to M^*$ 
such that $d_0 = e^*\zeta c.$ It follows that $d_0^* = c^*\zeta^* e^{**}.$
We have the commutative diagram
$$ 
\CD
  Q_{-1}^* @>e>> M \cr
  @V\phi=1 VV  @VV\phi_M V \cr
  Q_{-1}^* @>e^{**}>> M^{**}
\endCD
$$
(with $\phi = \phi_{Q_{-1}^*}$ the identity map), thus $e^{**} = \phi_Me$ and
therefore $d_0^* = c^*\zeta^* e^{**} = c^*\zeta \phi_Me.$

(ii) implies (i). We assume that $\zeta\:N\to M^*$ is an isomorphism with 
$d_0^* = c^*\zeta^* \phi_Me = (\zeta c)^*\phi_Me.$ 
Then $d_0 = d_0^{**} = e^*(\phi_M)^*(\zeta c)^{**}.$ There is the commutative diagram
$$
\CD
 M^* @<\zeta c<< Q_0 \cr
  @V\phi_{M^*}VV  @VV\phi = 1V \cr
 M^{***} @<(\zeta c)^{**}<< Q_0^{**},
\endCD
$$
therefore $(\zeta c)^{**} = \phi_{M^*}\zeta c$, thus 
$d_0 = e^*(\phi_M)^*(\zeta c)^{**} = e^*(\phi_M)^*\phi_{M^*}\zeta c =
e^*\zeta c$.
(Here we use that $(\phi_M)^*\phi_{M^*}$ is the identity map of $M^*.$ It implies
that $\phi_{M^*}$ is a splitting monomorphism; but in general,  $\phi_{M^*}$ is not
an isomorphism.)
This shows that $d_0$ is the composition of the cokernel map $\zeta c$ for $d_1$ with
the kernel map $e^*$ for $d^*_{-1}$. It follows that the sequence $Q_\bullet$ is exact. 
$\s$
	\medskip
{\bf 2.3.} Proof of Main Theorem.
Assume that $M$ is a reduced module 
such that both $M$ and $M^*$ are semi-Gorenstein-projective.
Let 
$$ 
 \cdots \to P_1 \to P_0 @>e>> M \to 0
$$ 
be a minimal projective
resolution. Since $M$ is semi-Gorenstein-projective, the $A$-dual sequence 
$$
 \cdots @<<< P_1^* @<<< P_0^* @<e^*<< M^{*} @<<< 0 \quad
$$ 
is exact.
Let 
$$
 \quad  0 @<<< M^* @<q<< Q_0 @<<< Q_1 @<<< \cdots
$$ 
be a minimal projective resolution 
of the right module $M^*$. The concatenation is an acyclic minimal 
complex $Q_\bullet$
of projective right modules (with $Q_{-i} = P_{i-1}^*$ for $i\ge 1$):
$$
 \cdots @<<< P_1^* @<<< P_0^* @<e^*q<< Q_0 @<<< Q_1 @<<<  \cdots,
$$
 Let us consider the $A$-dual $P_\bullet = Q_\bullet^*$
$$
 \cdots @>>> P_1 @>>> P_0 @>q^*e^{**}>> Q_0^* @>>> Q_1^* @>>> \cdots. \tag{$*$}
$$
It is the concatenation of the sequence
$$
 \cdots @>>> P_1 @>>> P_0 @>e^{**}>> M^{**} @>>> 0
$$
with the exact sequence 
$$
 0 @>>> M^{**} @>q^{*}>> Q_0^{*}  @>>> Q_1^* @>>> \cdots.
$$
In particular, the complex $(*)$ 
is exact at the positions  $P_i$  and $Q_i^*$ with $i\ge 1.$ 
	\medskip 
Conversely, let $P_\bullet = (P_i,f_i)_i$ be a minimal complex of projective modules
$$ 
 \cdots @>>> P_2 @>f_2>> P_1 @>f_1>> P_0 @>f_0>> P_{-1} @>f_{-1}>> P_{-1} @>>> \cdots 
$$ 
such that $H_i(P_\bullet) = 0$ for all $i\neq 0,-1,$ and such that the $A$-dual complex
$P_\bullet^*$ is acyclic. Let $e\:P_0 \to M$ be the cokernel of
$f_1$, thus 
$$ 
 \cdots @>>> P_2 @>f_2>> P_1 @>f_1>> P_0 @>e>> M \to 0
$$ 
is a minimal projective resolution of $M$. Since the complex $P_\bullet^*$ is acyclic,
it follows that $M$ is semi-Gorenstein-projective. 

Since $M$ is the cokernel
of $f_1$, we see that $M^*$ is the kernel of $f_{-1}^*.$
Since the complex $P_\bullet^*$ is acyclic, there is the exact sequence
$$
   0 @<<< M^* @<<< P_{-1}^* @<f_{-1}^*<< P_{-2}^* @<f_{-2}^*<< P_{-3}^* @<<< \cdots,
$$ 
and this is a projective resolution of $M^*$. Since the $A$-dual sequence 
$$
   P_{-1} @>f_{-1}>> P_{-2} @>f_{-2}>> P_{-3} @>>> \cdots.
$$ 
is exact, we see that $M^*$ is semi-Gorenstein-projective. 
$\s$
	\bigskip
{\bf 2.4. The $\infty$-torsionfree right modules $Z$ with $\Omega^2 Z$ 
semi-Gorenstein-projective.}	
We recall that a module $M$ is said to be {\it $\infty$-torsionfree} provided $\Tr M$
is semi-Gorenstein-projective. 
	\medskip
{\bf Proposition.} {\it The transpose $\Tr$ provides a bijection between the reduced
modules $M$ such that both $M$ and $M^*$ are semi-Gorenstein-projective and
the reduced $\infty$-torsionfree right modules $Z$ 
with $\Omega^2Z$ semi-Gorenstein-projective.}
	\medskip
For the proof, we need the following (well-known) lemma.
	\medskip
{\bf Lemma.} {\it For any module $M$, we have $\Omega^2 M = (\Tr M)^*$.}
	\medskip
Proof of Lemma. 
Take a minimal projective presentation 
$ P_1 @>f_1>> P_0 @>>> M  @>>> 0$
Then $\Omega^2 M  = \Ker f_1$.
By definition of $\Tr M$, there is the exact
sequence $P^*_0 @>f^*_1>> P^*_1 @>>> \Tr M \to 0$. If we apply $^* = \Hom(-,A_A)$, 
we get the exact sequence $0 \to (\Tr M)^*  @>>> P_1^{**} @>f_1^{**}>> P_0^{**}.$  
But $f_1^{**}$ can be identified with $f_1$, thus  
$(\Tr M)^* = \Ker f_1^{**} = \Ker f_1 = \Omega^2 M.$
$\s$
	\medskip
Proof of Proposition. 
Let $M$ be a reduced module such that both 
$M$ and $M^*$ are semi-Gorenstein-projective and $Z = \Tr M.$
Since $M = \Tr Z$ is semi-Gorenstein-projective, $Z$ is $\infty$-torsionfree. 
According to the lemma, $\Omega^2 Z  = (\Tr Z)^* = (\Tr\Tr M)^* = M^*$,
thus $\Omega^2 Z$ is semi-Gorenstein-projective.

Conversely, let $Z$ be a reduced $\infty$-torsionfree module such that $\Omega^2 Z$
is semi-Gorenstein-projective. Then $M = \Tr Z$ is semi-Gorenste-projective. 
The Lemma asserts that $M^* = (\Tr\Tr M)^* = (\Tr Z)^* = \Omega^2 Z$. This shows that
$M^*$ is semi-Gorenstein-projective.
$\s$
	\bigskip
{\bf 2.5 Summary.}
Let $M$ be a reduced module such that both 
$M$ and $M^*$ are semi-Gorenstein-projective.
The Main Theorem  yields a minimal complex $P_\bullet$.
Let us display, first, the complex $P_\bullet$ indicating the homology groups above,
and second, directly below, the acyclic $A$-dual complex $P_\bullet^*$.
We insert the modules $M$, $M^*,$ $M^{**}$, together
with the canonical map $\phi_M\:M \to M^{**}$ (shown as a bold
arrow), as well as the modules $\Tr M$ and $\Tr M^*$.
Since the modules $M$ and $M^*$ both are semi-Gorenstein-projective, 
the modules $\Tr M$ and $\Tr M^*$ are $\infty$-torsionfree. 
The complexes $P_\bullet$ and $P_\bullet^*$ provide minimal projective resolutions 
of the modules $M$ and $M^*$, respectively 
(they are encompassed by solid lines, with label sGp added). 
Similarly, $P_\bullet$ and $P_\bullet^*$ provide minimal projective coresolutions 
of the modules $\Tr M^*$ and $\Tr M$  which are concatenations of
$\mho$-sequences, respectively (these coresolutions are encompassed by dashed lines, with
label $\infty$-tf added).

$$
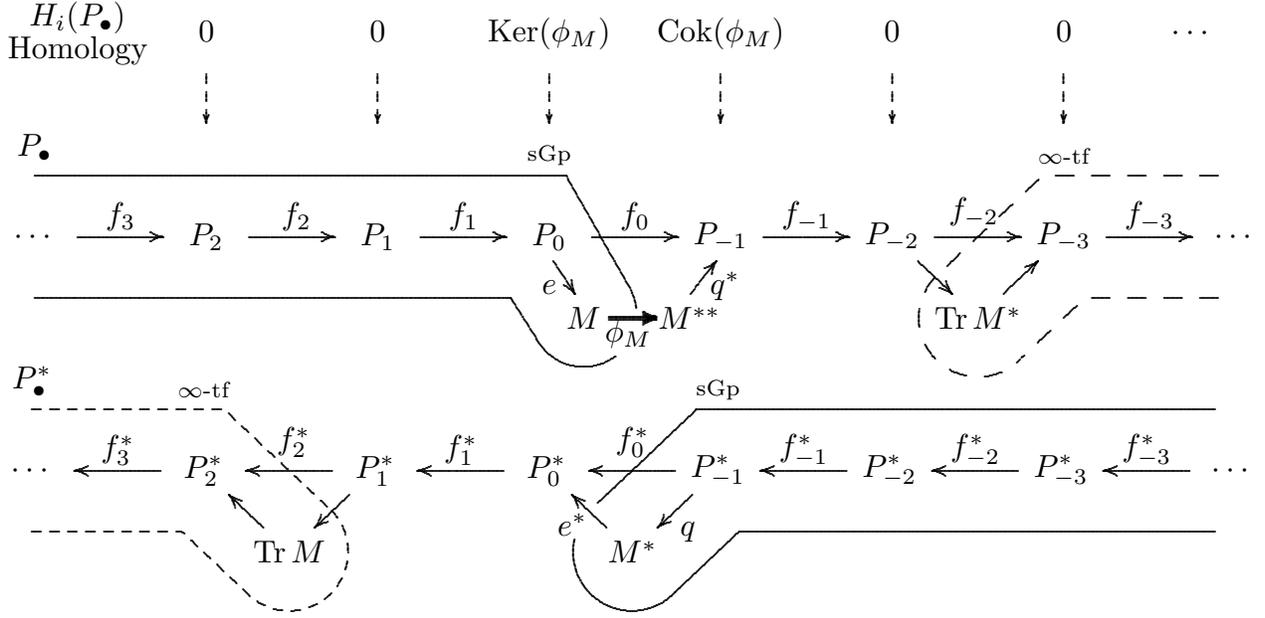
{\beginpicture
    \setcoordinatesystem units <.95cm,.9cm>

\put{\beginpicture
\put{$P_2$} at -5 0
\put{$P_1$} at -3 0
\put{$P_0$} at -1 0
\put{$P_{-1}$} at 1 0
\put{$P_{-2}$} at 3 0
\put{$P_{-3}$} at 5 0

\put{$M$} at -0.6 -1
\put{$M^{**}$} at 0.63 -1
\put{$e$\strut} at -1 -0.6
\put{$\phi_M$} at -.08 -1.2
\put{$q^*$\strut} at 1.05 -0.6

\plot -7 .75  -.8 .75  -.05 -.65 /
\plot -7 -.75 -1.45 -.75 -1.15 -1.25 /
\circulararc 100 degrees from -1.15 -1.25 center at -.6 -.95 
\circulararc -20 degrees from -.04 -.67 center at -.6 -.95 
\put{$\ss \text{sGp}$} at -1 1

\put{$\cdots$} at -7 0
\put{$\cdots$} at 7 0
\arr{-6.5 0}{-5.5 0}
\arr{-4.5 0}{-3.5 0}
\arr{-2.5 0}{-1.5 0}
\arr{-.5 0}{.5 0}
\arr{1.5 0}{2.5 0}
\arr{3.5 0}{4.5 0}
\arr{5.5 0}{6.5 0}

\arr{-0.95 -.3}{-.7 -.7}
\arr{-0.3 -1}{.25 -1}
\arr{-0.3 -1.02}{.25 -1.02}
\arr{-0.3 -.98}{.25 -.98}

\arr{.65 -.7}{0.95 -.3}

\put{$f_3$} at -6 0.25
\put{$f_2$} at -3.95 0.25
\put{$f_1$} at -2 0.25
\put{$f_0$} at 0 0.25
\put{$f_{-1}$} at 2 0.25
\put{$f_{-2}$} at 3.95 0.3
\put{$f_{-3}$} at 6 0.25

\setdashes <2mm>
\plot 6.8 .75   4.75 .75   3.5 -.5 /
\plot 6.8 -.75  5.25 -.75  4.5 -1.5 /
\circulararc 180 degrees from 3.5 -.5 center at 4 -1 
\put{$\ss \text{$\infty$-tf}$} at 5 1
\setsolid

\put{$\Tr M^*$} at 4 -1 
\arr{3.3 -0.3}{3.7 -.7}
\arr{4.3 -0.7}{4.7 -.3}

\setdashes <1mm>
\arr{-5 2.}{-5 1.4}
\arr{-3 2.}{-3 1.4}
\arr{-1 2.}{-1 1.4}
\arr{1 2.}{1 1.4}
\arr{3 2.}{3 1.4}
\arr{5 2.}{5 1.4}
\put{Homology} at -6.5 2.3
\put{$H_i(P_\bullet)$} at -6.5 2.7
\multiput{$0$\strut} at -5 2.5  -3 2.5  3 2.5  5 2.5 /
\put{$\Ker(\phi_M)$\strut} at -1 2.5
\put{$\Cok(\phi_M)$\strut} at 1 2.5
\put{$\cdots$} at 6.5 2.5 

\put{$P_\bullet$} at -7 1.1

\endpicture} at 0 0

\put{\beginpicture
\put{$P_2^*$} at -5 0
\put{$P_1^*$} at -3 0
\put{$P_0^*$} at -1 0
\put{$M^*$} at 0 -1
\put{$P_{-1}^*$} at 1 0
\put{$P_{-2}^*$} at 3 0
\put{$P_{-3}^*$} at 5 0
\put{$q$\strut} at 0.65 -.7
\put{$e^*$\strut} at -.7 -.7

\put{$\Tr M$} at -4 -1

\setdashes <1mm>
\plot -7 .75  -4.75 .75  -3.5 -.5 /
\plot -7 -.75 -5.25 -.75 -4.5 -1.5 /
\circulararc -180 degrees from -3.5 -.5 center at -4 -1 
\put{$\ss \text{$\infty$-tf}$} at -5 1
\setsolid

\plot 6.8 .75   0.75 .75   -.4 -.4 /
\plot 6.8 -.75  1.25 -.75  .5 -1.5 /
\circulararc -149 degrees from .5 -1.5 center at 0 -1 
\put{$\ss \text{sGp}$} at 1 1

\put{$\cdots$} at -7 0
\put{$\cdots$} at 7 0
\arr{-5.5 0}{-6.5 0}
\arr{-3.5 0}{-4.5 0}
\arr{-1.5 0}{-2.5 0}
\arr{.5 0}{-.5 0}
\arr{2.5 0}{1.5 0}
\arr{4.5 0}{3.5 0}
\arr{6.5 0}{5.5 0}

\arr{-4.3 -0.7}{-4.7 -.3}
\arr{-3.3 -0.3}{-3.7 -.7}


\arr{-0.3 -0.7}{-0.7 -.3}
\arr{.7 -0.3}{.3 -.7}
\put{$f_3^*$} at -6 0.25
\put{$f_2^*$} at -3.95 0.35
\put{$f_1^*$} at -2 0.25
\put{$f_0^*$} at 0 0.35
\put{$f_{-1}^*$} at 2 0.25
\put{$f_{-2}^*$} at 4 0.25
\put{$f_{-3}^*$} at 6 0.25

\put{$P_\bullet^*$} at -7 1.1
\put{} at 0 -1.4
\endpicture} at 0 -3.7

\endpicture}
$$
	\bigskip
Be aware that the complexes $P_\bullet$ and $P_\bullet^*$ with the 
accompagnying modules seem to look quite similar, however there is a decisive
difference: whereas the complex $P_\bullet^*$ is acyclic, the complex $P_\bullet$ 
usually is not acyclic (its homology modules are mentioned above the complex). 
Let us stress that $P_\bullet$ is acyclic if and only if $M$ is Gorenstein-projective. 
	\medskip
{\bf 2.6. Remarks. (1)} The Main Theorem illustrates nicely that {\it an indecomposable module $M$ is 
Gorenstein-projective if and only if both $M$ and $M^*$ are semi-Gorenstein-projective
and $M$ is reflexive} (since the latter means that $\phi_M$ is an isomorphism), as known from [AB],
and stressed for example in [AM].
	\medskip 
{\bf (2)} By construction, the complex $P_\bullet$ (and thus also $P_\bullet^*$)
is uniquely determined by the module $M$. Let us stress that {\it $P_\bullet$ is usually 
{\bf not} determined by $M^*$}. 

In general, given an acyclic minimal complex 
$Q_\bullet = (Q_i,d_i\:Q_i\to Q_{i-1})$ of projective right modules (such as 
$Q_\bullet = P_\bullet^*$),
say with $N_i$ being the image of $d_i$, then any module $N_i$ determines uniquely
the modules $N_{j}$ with $j\le i$, since $N_{j} = \Omega^{j-i}N_i,$ but usually not the
modules $N_j$ with $j > i.$ 

If we look at the complexes $P_\bullet$ which are obtained in Main Theorem,
then there do exist examples, where $P_\bullet = (P_i,f_i)$ 
is {\bf not} determined by $M^*$ (this is the image of $f_0$), as shown in [RZ2]. 
Namely, let $q \in k$ be an element with infinite 
multiplicative order and $A = \Lambda(q)$ the algebra defined in
[RZ2](1.1). Then the modules $M$ of the form $M(1,-q,c)$ with $c\in k$ are indecomposable,
non-projective and semi-Gorenstein-projective and they are pairwise non-isomorphic
(thus also the right modules $\Tr M = \Tr M(1,-q,c)$ are pairwise non-isomorphic) ---
whereas all the right modules $M^* = M(1,-q,c)^*$ are isomorphic, and also
semi-Gorenstein-projective, see [RZ2](1.7); they are of the form $M^* = M'(1,-q^{-1},0)$,
see [RZ2](9.4). Actually, in this case already all the right
modules $\Omega \Tr M = (\Omega M)^*$ are isomorphic, namely of the form $M'(1,-1,0)$, see 
[RZ2](3.2). 

To phrase it differently: [RZ2] provides an infinite family of acyclic minimal complexes
$Q(c)_\bullet = (Q(c)_i,d(c)_i)$ 
indexed by the elements $c\in k$, such that for any $i\in \Bbb Z$, 
the images of the maps $d(c)_i$ are pairwise non-isomorphic if $i\ge 0$,
but pairwise isomorphic if $i < 0$. 

	\medskip
{\bf (3)} Even if the modules $M$ and $M^*$ are indecomposable, the module $M^{**}$
may be decomposable, as the example of $M = M(q)$ in [RZ1] shows. Note that if $M^*$ 
is indecomposable and not projective (this is the case in the example), then also
$\Tr M^*$ is indecomposable and not projective, thus in the complex $P_\bullet$ displayed 
in 2.5, the images of all the maps $f_i$ with $i\neq 0, -1$ can be 
indecomposable and not projective, whereas $M^{**}$ is decomposable. 
	\medskip
{\bf  (4)} Let $A$ be a connected algebra with a
non-reflexive module $M$ such that both $M$ and $M^*$ are semi-Gorenstein-projective.
Then, of course $A$ is not left weakly Gorenstein (recall that an algebra $A$
is said to be {\it left weakly Gorenstein,} provided any semi-Gorenstein-projective
module is Gorenstein-projective, see [RZ1]). Is it possible that $A$ is right weakly
Gorenstein (this means that any $\infty$-torsionfree module is Gorenstein-projective)?
Of course, as we have mentioned in 2.4, the module $\Tr M^*$ is always 
$\infty$-torsionfree. Thus, if $\Tr M^*$ is not Gorenstein-projective, then 
$A$ is not right weakly Gorenstein. But we do not know whether $\Tr M^*$ can
be Gorenstein-projective. 

In section 3, we discuss the extreme case that $M^*$ is projective
(thus $\Tr M^* = 0$). According to the classical homological conjectures,
this case should be impossible, see 3.5. But could $M^*$ be Gorenstein-projective?
	\medskip
{\bf (5)} After completing the paper, the authors became aware of the recent preprint [G] by
G\'elinas which also deals with the complex $P_\bullet$. There, the central
(and decisive) map $q^*\phi_Me$
is called the {\it Norm map} of the module $M$, with reference to Buchweitz [B], 5.6.1. 
	\bigskip\bigskip
{\bf 3. Special cases.}
	\medskip 
{\bf 3.1. The case where $\phi_M$ is an epimorphism (or a monomorphism).}
Let us consider now the special case of a module $M$ with both $M,M^*$
semi-Gorenstein-projective such that $\phi_M$ is an epimorphism (or a
monomorphism). But
we stress from the beginning that at present {\bf no non-trivial such example is
known} (all known modules $M$ with both $M,M^*$ semi-Gorenstein-projective 
such that $\phi_M$ is an epimorphism or a
monomorphism, are Gorenstein-projective).
	\bigskip
{\bf Proposition.} 
	\smallskip
(1) {\it Let $M$ be a semi-Gorenstein-projective module. Then we have:
$M^*$ is semi-Gorenstein-projective and $\phi_M$ is an epimorphism if and only if 
$(\Omega M)^*$ is semi-Gorenstein-projective.}
	
(2) {\it Let $M'$ be a torsionless semi-Gorenstein-projective module and let $M = \mho M'.$
Then $M$ is semi-Gorenstein-projective. And $(M')^*$ is semi-Gorenstein-projective if 
and only if $M^*$ is also semi-Gorenstein-projective and $\phi_M$ is an epimorphism.}
	\medskip
Here, we consider in (1) a module $M$ such that both $M$ and $M^*$ are semi-Gorenstein-projective,
and in (2), a module $M'$ such that both $M'$ and $(M')^*$ are
semi-Gorenstein-projective. In (1) we deal with the case that $\phi_M$ is an epimorphism. 
In (2) we deal with the case that $\phi_{M'}$ is an monomorphism (namely,
$M'$ is torsionless if and only if $\phi_{M'}$ is a monomorphisms).
	\bigskip
Proof of (1). We can assume that $M$ is indecomposable and not projective. 
Since $M$ is semi-Gorenstein-projective, it follows that the exact sequence
$0 \to \Omega M \to P(M) \to M \to 0$ is an $\mho$-sequence, and
$0 \to M^* \to P(M)^* \to (\Omega M)^* \to 0$ is exact. Thus, $(\Omega M)^*$ is 
semi-Gorenstein-projective if and only if $M^*$ is semi-Gorenstein-projective and
$\Ext^1((\Omega M)^*,A_A) = 0.$ According to Lemma 2.4(b) in [RZ1], we have 
$\Ext^1((\Omega M)^*,A_A) = 0$ if and only if $\phi_M$ is an epimorphism. 
	\smallskip
Proof of (2). We can assume that $M'$ is indecomposable and not projective. 
Since $M'$ is torsionless and semi-Gorenstein-projective, the module $M = \mho M'$
is semi-Gorenstein-projective. There is an $\mho$-sequence $0 \to M' \to P \to M \to 0$,
thus an exact sequence $0 \to M^* \to P^* \to (M')^* \to 0.$ Then $(M')^*$ is
semi-Gorenstein-projective if and only if $M^*$ is semi-Gorenstein-projective and
$\Ext^1((M')^*,A_A) = 0.$ According to 2.4(b) of [RZ1], we have
$\Ext^1((M')^*,A_A) = 0$ if and only if $\phi_M$ is an epimorphism.
$\s$
	\medskip
{\bf Remark.} If we denote by $\Cal M$ the class of reduced
modules $M$ such that both $M$ and $M^*$ are semi-Gorenstein-projective and $\phi_M$ is
an epimorphism, and by $\Cal M'$ the class of reduced
modules $M$ such that both $M'$ and $(M')^*$ are semi-Gorenstein-projective and $\phi_{M'}$ is
a monomorphism, then $\Omega$ and $\mho$ provide
inverse bijections between isomorphism classes as follows:
$$
{\beginpicture
\setcoordinatesystem units <2cm,1cm>
\put{$\Cal M'$} at 0 0
\put{$\Cal M$} at 2 0
\arr{0.5 0.1}{1.5 0.1}
\arr{1.5 -.1}{0.5 -.1}
\put{$\Omega$} at 1 -.3
\put{$\mho$} at 1 .3
\endpicture}
$$
If $M$ belongs to $\Cal M$ and $M' = \Omega M$ (thus $\mho M' = M$), then 
$\Cok \phi_{M*} \simeq \Ker \phi_M$. 
	\bigskip
{\bf 3.2. The semi-Gorenstein-projective modules $M$ with $M^*$ projective.}
Another special case should be considered, namely the case of a 
semi-Gorenstein-projective reduced module such that 
$M^*$ is projective. Also here, let us stress 
from the beginning that at present {\bf no non-trivial such example is
known} (all known semi-Gorenstein-projective modules $M$ with $M^*$
projective are Gorenstein-projective, thus even projective).

 Let $M$ be a semi-Gorenstein-projective module with $M^*$ being projective.
In addition, we may assume that $M$ is reduced.
Since $M^*$ is projective, we take as projective cover $q\:Q_0 \to M^*$ the 
identity map $1 = 1_{M^*}$; thus $f_0 = \phi_M\cdot e.$
The diagram considered in 2.5 now has the following
special form:
\vfill\eject
$$
{\beginpicture
    \setcoordinatesystem units <.95cm,.9cm>

\put{\beginpicture
\put{$P_2$} at -5 0
\put{$P_1$} at -3 0
\put{$P_0$} at -1 0
\put{$M^{**}$} at 1 0
\put{$0$} at 3 0
\put{$0$} at 5 0

\put{$M$} at -0 -1
\put{$e$\strut} at -.75 -0.72
\put{$\phi_M$\strut} at .85 -.75


\put{$\cdots$} at -7 0
\put{$\cdots$} at 7 0
\arr{-6.5 0}{-5.5 0}
\arr{-4.5 0}{-3.5 0}
\arr{-2.5 0}{-1.5 0}
\arr{-.5 0}{.5 0}
\arr{1.5 0}{2.5 0}
\arr{3.5 0}{4.5 0}
\arr{5.5 0}{6.5 0}

\arr{-0.8 -.3}{-.25 -.85}
\arr{.3 -.8}{0.8 -.3}

\put{$f_3$} at -6 0.25
\put{$f_2$} at -3.95 0.25
\put{$f_1$} at -2 0.25
\put{$f_0$} at 0 0.25
\put{$$} at 2 0.25
\put{$$} at 3.95 0.3
\put{$$} at 6 0.25

\plot -7 .75  -.85 .75  .4 -.5 /
\plot -7 -.75 -1.25 -.75 -.5 -1.5 /
\circulararc  140 degrees from -.5 -1.5 center at 0 -1 
\put{$\ss \text{sGp}$} at -1.4 -1.2
\setsolid

\setdashes <1mm>
\arr{-5 2.}{-5 1.4}
\arr{-3 2.}{-3 1.4}
\arr{-1 2.}{-1 1.4}
\arr{1 2.}{1 1.4}
\arr{3 2.}{3 1.4}
\arr{5 2.}{5 1.4}
\put{Homology} at -6.5 2.3
\put{$H_i(P_\bullet)$} at -6.5 2.7
\multiput{$0$\strut} at -5 2.5  -3 2.5  3 2.5  5 2.5 /
\put{$\Ker(\phi_M)$\strut} at -1 2.5
\put{$\Cok(\phi_M)$\strut} at 1 2.5
\put{$\cdots$} at 6.5 2.5 

\put{$P_\bullet$} at -7 1.1

\endpicture} at 0 -4.2

\put{\beginpicture
\put{$P_2^*$} at -5 0
\put{$P_1^*$} at -3 0
\put{$P_0^*$} at -1 0
\put{$M^*$} at 1 0
\put{$0$} at 3 0
\put{$0$} at 5 0

\put{$\Tr M$} at -4 -1

\setdashes <1mm>
\plot -7 .75  -4.75 .75  -3.5 -.5 /
\plot -7 -.75 -5.25 -.75 -4.5 -1.5 /
\circulararc -180 degrees from -3.5 -.5 center at -4 -1 
\put{$\ss \text{$\infty$-tf}$} at -5 1
\setsolid

\put{$\cdots$} at -7 0
\put{$\cdots$} at 7 0
\arr{-5.5 0}{-6.5 0}
\arr{-3.5 0}{-4.5 0}
\arr{-1.5 0}{-2.5 0}
\arr{.5 0}{-.5 0}
\arr{2.5 0}{1.5 0}
\arr{4.5 0}{3.5 0}
\arr{6.5 0}{5.5 0}

\arr{-4.3 -0.7}{-4.7 -.3}
\arr{-3.3 -0.3}{-3.7 -.7}


\put{$f_3^*$} at -6 0.25
\put{$f_2^*$} at -3.95 0.35
\put{$f_1^*$} at -2 0.25
\put{$e^*$} at 0.1 0.3
\put{$$} at 2 0.25
\put{$$} at 4 0.25
\put{$$} at 6 0.25
\put{$P_\bullet^*$} at -7 1.1
\endpicture} at 0 -8

\endpicture}
$$
	\bigskip
\noindent
Here, $P_\bullet^*$ is acyclic, $\Tr M$ is the image of $f_{-2}^*$ and 
$e\:M^* = P_0^*$ is an inclusion map.
It follows that $\Tr M$ has projective dimension at most 2
(and projective dimension at most 1, in case $M^* = 0$). 
	\medskip
{\bf 3.3.} Recall that a module $M$ is a Nunke module 
provided $M$ is semi-Gorenstein-projective
and $M^* = 0$ (and the Nunke condition for an algebra $A$ asserts that
the zero module is the only Nunke $A$-module, see 3.5). 
	\medskip
{\bf Proposition.} {\it The transpose $\Tr$ provides a bijection between 
the modules $M$ which are semi-Gorenstein-projective, with $M^*$ being projective,
and the $\infty$-torsionfree right modules $Z$ of projective dimension at most $2$.

The transpose $\Tr$ provides a bijection between 
the Nunke modules $M$ and the $\infty$-torsionfree 
right reduced modules $Z$ of projective dimension at most $1$.}
$\s$
	\medskip
{\bf 3.4. Nunke modules.}
We consider now Nunke modules. As we will see,  
this is essentially just the situation
where, on the one hand, $\phi_M$ is an epimorphism (as discussed in 3.1), and,
on the other hand, $M^*$ is projective (as discussed in 3.2). Of course, again we have to stress
that {\bf no non-trivial example is known} (the only known Nunke module is the zero module)
and there is the old conjecture (one of the classical homological conjectures) that 
no non-zero Nunke module exists (for a further discussion of corresponding
conjectures, see 3.5).

	\medskip
{\bf Proposition.} {\it Let $M$ be a module. The following conditions are equivalent:
\item{\rm (i)} $M$ is semi-Gorenstein-projective, 
  $\phi_M$ is an epimorphism and $M^*$ is projective.
\item{\rm (ii)} $M$ is the direct sum of a projective module and a Nunke module.}
	\medskip
Proof. We can assume that $M$ is indecomposable. 

(i) $\implies$ (ii). 
Let $M$ be a semi-Gorenstein-projective module, with $\phi_M$ an epimorphism
and $M^*$ projective.
Then also $M^{**}$ is projective, thus $\phi_M$ is a split epimorphism, thus 
$M \simeq M'\oplus M^{**}$. Since we assume that
$M$ is indecomposable, there are two possibilities: 
Either $M' = 0$, thus $M \simeq M^{**}$, thus $M$ is projective. 
Or else $M^{**} = 0$ and thus already $M^* = 0$
(namely, if $M^* \neq 0$, then also $M^{**} \neq 0,$ since the module
$M^*$ is torsionless), thus $M$ is a Nunke module.

(ii) $\implies$ (i). If $M$ is projective, then $M$ is
semi-Gorenstein-projective. Also, $M$ is reflexive, thus  $\phi_M$ is surjective, and
with $M$ projective, also $M^*$ is projective.. If $M$ is a 
Nunke module, then $M^* = 0$ implies that $M^{**} = 0$, thus $\phi_M$ is surjective and
$M^*$, being the zero module, is projective. 
$\s$
	\bigskip
{\bf 3.5. Some conjectures.}
	\medskip 
{\bf Proposition.} 
{\it Let $A$ be an artin algebra. We consider the following
conditions:
	\smallskip 
\item{\rm(1)} There is a bound $b$ with the following property: If $M$ is a right
  $A$-module of finite projective dimension, then the projective dimension of $M$ is at most
  $b$.
\item{\rm(2)} A semi-Gorenstein-projective $A$-module $M$ with $M^*$ of finite projective dimension is projective.
\item{\rm(3)} A semi-Gorenstein-projective $A$-module $M$ with $M^*$ projective is projective.
\item{\rm(4)} The only Nunke module is the zero module
\item{\rm(5)} There is no simple module which is a Nunke module.
\item{\rm(6)} There is no simple injective module which is a Nunke module.
	\smallskip 
Then} (1) $\implies$ (2) $\implies$ (3) $\implies$ (4) $\implies$ (5) $\implies$ (6).
	\medskip
{\bf Remark.} (1) is called the {\bf finitistic dimension conjecture} for $A^{\op}$,
and (4) is called the {\bf Nunke condition} for $A$. Both are classical homological 
conjectures (that these conditions hold true for all finite-dimensional algebras)
and it is well-known that (1) implies (4). See, for example, [H].
The assertion (5) is called the {\bf weak Nunke condition}, it is equivalent to the
{\bf generalized Nakayama conjecture} (see [AR], Proposition 1.5).
The conditions (1), (4), and (5) are mentioned in [ARS] as
conjectures 11, 12, and 9, respectively). The Proposition formulates the intermediate conjectures
(2) and (3). Also, we mention the weaker conjecture (6) which will be discussed in
3.6 and 3.7.
	\medskip
Proof of Proposition. (1) $\implies$ (2). 
Let $M$ be semi-Gorenstein-projective.
We may assume that $M$ is indecomposable and not projective. Now
$Z = \Tr M$ is $\infty$-torsionfree, thus there is an
exact sequence 
$$
 0 \to Z \to P^0 @>d^0>>P^1 @>d^1>> \cdots
$$
with projective modules $P^i$, where $i\ge 0.$ 
Let $M^i$ be the image of $d^{i-1}$ for all $i\ge 1.$ Since $Z$ is 
indecomposable and not projective, we see that $\Omega^i M^i = Z$ 
for all $i\ge 1.$ 

We apply the Lemma from 2.4 to $Z = \Tr M$ and see that $\Omega^2Z = 
(\Tr\Tr M)^* = M^*,$ since $M$ is reduced. 

By assumption, $M^*$ has finite projective dimension, thus $\Omega^a (M^*) = 0$, for
some $a\ge 0$ and therefore $\Omega^{i+2+a}M^i =  \Omega^{2+a} Z = \Omega^a M^* = 0.$
Thus $M^i$ has finite projective dimension, for all $i\ge 1$. 
According to (1), we know that $M^i$ has projective dimension at most $b$, for all $i\ge 1.$
Since $M^{b+1}$ has projective dimension at most $b$, we see that $Z = \Omega^{b+1}M^{b+1} = 0.$
But $Z = 0$ imples that $M = \Tr Z = 0.$
	\medskip
(2) $\implies$ (3) is trivial. (3) $\implies$ (4): Let $M$ be semi-Gorenstein-projective
and $M^* = 0.$ Since $M^*$ is projective, it follows from (3) that $M$ is projective.
But a projective module $P$ with $P^* = 0$ is the zero module. The implications
(4) $\implies$ (5)  $\implies$ (6) are again trivial. 
$\s$
	\medskip
{\bf Remark.} As we have mentioned already, the weak Nunke condition 3.5 (5)
is equivalent to the following conjecture, now called the {\bf Auslander-Reiten conjecture:}
{\it There is no non-zero semi-Gorenstein-projective $M$, with $\Ext_A^i(M,M) = 0$ for all $i\ge 1.$}
This was shown by Auslander and Reiten [AR]. 

The weak Nunke condition asserts that a simple semi-Gorenstein-projective module should be 
torsionless. One may ask whether a simple semi-Gorenstein-projective module should be even
Gorenstein-projective, but at present, nothing is known about simple semi-Gorenstein-projective
modules, this is really a pity. There is the special case of a simple injective module. Is it
possible that a simple injective module is semi-Gorensetin-projective without being 
already projective (thus uninteresting)? In 3.5, we have added the corresponding conjecture: 
{\it There are no simple injective Nunke modules}, as conjecture (6). 
One may call this conjecture (6) the 
{\bf very weak Nunke condition.} In 3.6 and 3.7, we
show that the very weak Nunke condition is equivalent to 
a weak form of the Auslander-Reiten conjecture.
	\bigskip
{\bf 3.6. Simple injective semi-Gorenstein-projective modules.}
The conjecture 3.5 (6) asserts that a simple injective semi-Gorenstein-projective module $S$
should be projective (thus $\Omega S = 0$). In 3.6, we look at a simple injective
semi-Gorenstein-projective module and try to analyse $\Omega S$.
Let $A$ be an artin algebra and $S$ a simple injective module with
endomorphism ring $\End(S) = D.$
Let $M = \Omega S$. Let $e$ be a primitive idempotent
of $A$ such that $eS\neq 0.$ Let $B = A/AeA$, thus $\Cal U(S) = \mod B$ is a full subcategory of 
$\mod A$; it consists of all the $A$-modules which do not have $S$ as a composition factor. 
Note that $M$ belongs to $\Cal U(S)$. 

As we have mentioned, $A$ is supposed to be an artin algebra. 
We assume that $A$ is a $k$-algebra, where $k$ is a commutative artinian ring and 
$A$ is of finite length as a $k$-module. 
	\medskip
{\bf Proposition.} {\it Let $S$ be a simple injective module with $\End(S) = D.$
Let $M = \Omega S$. Let $e$ be a primitive idempotent of $A$ such that $eS\neq 0$ and
$B = A/AeA$. The following conditions are equivalent:
\item{\rm (i)} The module $S$ is semi-Gorenstein-projective.
\item{\rm (ii)} The $B$-module $M$ is a Nunke $B$-module, $\Ext_B^i(M,M) = 0$ for all $i\ge 1$,
  and either $M = 0$ or else $\End (M)$ is isomorphic to $D$ as a $k$-algebra.\par}
	\medskip 
We need some preparations for the proof.
	\medskip
{\bf (a)} {\it Projective $B$-module are projective when considered as $A$-modules.
Thus $\Cal U(S)$ is closed under projective covers.}  $\s$
	\medskip
{\bf (b)} {\it If $U,U'$ are $B$-modules, then $\Ext_A^i(U,U') = \Ext_B^i(U,U')$ for all $i\ge 0$.}
	\medskip
Proof. This is clear for $i = 0.$ For $i\ge 1$, 
we start with a projective resolution $P_\bullet$ of the $B$-module $U$.
According to (a), this is also a projective resolution of $U$ considered
as an $A$-module. We form the complex
$\Hom_B(P_\bullet,U') = \Hom_A(P_\bullet,U')$ and consider its homology.
$\s$
	\medskip
{\bf (c)} {\it If $U$ is a $B$-module, then $\Ext_A^i(U,P(S)) = \Ext_A^i(U,M) (= \Ext_B^i(U,M))$
for all $i\ge 0$.}
	\medskip
Proof. We apply $\Hom_A(U,-)$ to the exact sequence $0 \to M \to P(S) \to S \to 0.$ We have
$\Hom(U,S) = 0$, since $U$ belongs to $\Cal U(S)$, and we have $\Ext^i(U,S) = 0$ for $i\ge 1$,
since $S$ is injective. This shows that $\Ext_A^i(U,M) = \Ext_A^i(U,P(S))$ for all $i\ge 0.$
For the second equality, see (b).
$\s$
	\medskip
{\bf (d)} {\it Let $U$ be a $B$-module. Then $U$ is semi-Gorenstein-projective as
an $A$-module if and only if $U$ is semi-Gorenstein-projective as a $B$-module and
$\Ext^i_B(U,M) = 0$ for all $i\ge 1.$}
	\medskip
Proof. We can assume that $A$ is basic, thus  ${}_AA = B\oplus P(S)$. 
By definition, $U$ is semi-Gorenstein-projective as an $A$-module if and only if
$\Ext_A^i(U,B) = 0$ and $\Ext_A^i(U,P(S)) = 0$ for all $i\ge 1.$ According to (b), we have
$\Ext^i_B(U,B) = \Ext_A^i(U,B)$ for all $i\ge 1$. But $\Ext_A^i(U,B) = 0$ 
for all $i\ge 1$ means that $U$ is semi-Gorenstein-projective as a $B$-module. 
Also, by (c) we have  $\Ext_A^i(U,P(S)) =  \Ext_B^i(U,M)$ for all $i\ge 0$.
$\s$
	\medskip 
{\bf (e)} {\it Let $M\neq 0$. 
The embedding $M \to P(S)$ is a left $\add({}_AA)$-approximation
if and only if $\End (M)$ is isomorphic to $D$ as a $k$-algebra 
and $\Hom(M,{}_BB) = 0.$}
	\medskip
Proof. Let $\iota\:M \to P(S)$ be the inclusion map.
We assume again that $A$ is basic, thus  ${}_AA = B\oplus P(S)$. 

A map $f\:M \to {}_BB$ factors through $\iota$ if and only if
$f = 0$ (since $\Hom(P(S),{}_BB) = 0).$ Thus all maps 
$f\:M \to {}_BB$ factor through $\iota$ if and only if $\Hom(M,{}_BB) = 0.$

It remains to look at maps $M \to P(S)$. Any endomorphism of $P(S)$ maps $M = \rad P(S)$ into
$M$, thus there are canonical maps $\pi'\:\End(P(S)) \to \End(M)$ as well as
$\pi\:\End(P(S)) \to \End(S).$ Since $P(S)$ is a projective cover of $S$, thus
$\pi$ is surjective.
Since $S$ is injective, we have $\Hom(P(S),M) = 0$, $\pi$ is an injective map. Altogether, we
see that $\pi$ is an isomorphism. 
Since $M \neq 0$, and $S$ is injective, $\Hom(S,P(S)) = 0$, therefore $\pi'$ is injective.
It follows that $\mu = \pi'\pi^{-1}$ is  an embedding  of $D = \End(S) \simeq \End(P(S))$ 
into $\End(M).$ Of course, $\mu$ is an surjective if and only if $\iota$ is a
$\add({}_AA)$-approximation. Note that the embedding $\mu$ is an isomorphism if and only
if the length of ${}_kD$ is equal to the length of ${}_k\End(M)$, thus if and only if
$\End(M)$ is isomorphic to $D$ as a $k$-algebra.
$\s$
	\medskip
Proof of Proposition. (i) $\implies$ (ii). Let $S$ be semi-Gorenstein-projective.
If $S$ is projective, then $M = 0$ and the conditions in (ii) are trivially satisfied. 
Thus, we assume that $S$ is not projective, thus 
Then ${}_AM = \Omega S$ is an indecomposable 
semi-Gorenstein-projective $A$-module. According to (d),
$M$ considered as a $B$-module is semi-Gorenstein-projective and 
$\Ext^i_B(M,M) = 0$ for all $i\ge 1.$ Also, the inclusion $\iota\:M \to P(S)$ is a left
$\add({}_AA)$-approximation (since $\Ext^1(S,{}_AA) = 0$),
thus (e) asserts that 
$\End (M)$ is isomorphic to $D$ as a $k$-algebra
and $\Hom(M,{}_BB) = 0.$ In particular, we see that $M$ is a Nunke $B$-module.

(ii) $\implies$ (i). We can assume that $M$ is non-zero 
(otherwise $S$ is projective, thus of course
semi-Gorenstein-projective). We assume that $M$ is a Nunke
$B$-module, that $\End (M)$ is isomorphic to $D$ as a $k$-algebra
(in particular, $M$ is indecomposable) and that $\Ext_B^i(M,M) = 0$ for all $i\ge 1$.
Since $M$ is a semi-Gorenstein-projective
$B$-module and $\Ext_B^i(M,M) = 0$ for all $i\ge 1$, we 
can apply (d) to $U = M$ and see that $M$ is semi-Gorenstein-projective
also as an $A$-module.
Since $\End (M)$ is isomorphic to $D$ as a $k$-algebra 
and  $\Hom(M,{}_BB) = 0,$ we can use (e). It
asserts that the embedding $M \to P(S)$ is 
a left $\add({}_AA)$-approximation. Since $M$ is a Nunke $B$-module, we have $\Hom(M,{}_BB) = 0$.
It follows from $M\neq 0$ that $M$ is not a projective $B$-module. Of course, $M$ is also not
isomorphic to $P(S)$, thus $M$ is not projective as an $A$-module.
Since $M$ is semi-Gorenstein-projective,
indecomposable and not projective, it follows that $S = \mho M$ is semi-Gorenstein-projective.
$\s$
	\medskip
{\bf 3.7. Corollary.} {\it There exists an artin algebra $A$ with a simple injective Nunke module
if and only if there exists an artin algebra $B$ with an indecomposable 
semi-Gorenstein-projective module $M$ with $\rad\End(M) = 0$, 
such that $\Hom(M,{}_BB) = 0$ and $\Ext_B^i(M,M) = 0$ for all $i \ge 1.$}
	\medskip
Proof. Let $A$ be an artin algebra with a
simple, injective, semi-Gorenstein-projective module $S$ with $S^* = 0.$
Let 
$\Cal U(S) = \mod B$ for some artin algebra $B$. Then Proposition 3.6 asserts that 
the $B$-module $M$ is semi-Gorenstein-projective, 
$\Hom(M,{}_BB) = 0$, and $\Ext_B^i(M,M) = 0$ for all $i \ge 1,$
and $\End(M)$ is isomorphic to $\End(S)$ as a $k$-algebra. Of course, if $\End(M)$ is
isomorphic to $\End(S)$, then $\End(M)$ is a division ring, thus $M$ is indecomposable and
$\rad\End(M) = 0.$

Conversely, let $B$ be an artin algebra with an indecomposable $B$.module $M$ which is
semi-Gorenstein-projective, such that $\rad\End(M) = 0$, 
$\Hom(M,{}_BB) = 0$, and finally $\Ext_B^i(M,M) = 0$ for all $i \ge 1.$ 
Let $D = \End(M)$. Since $M$ is indecomposable, $D$ is a local artin algebra. Since
$\rad\End(M) = 0$, we see that $D$ is a division ring.
Let $A = \left[\smallmatrix B & M \cr 0 & D^{\op}\endsmallmatrix\right]$. Let $P$ be the
indecomposable projective $A$-module  $P = 
\left[\smallmatrix M \cr D^{\op}\endsmallmatrix\right]$, let $S = P/\rad P.$ Then $\End(S) = D$.
Note that $S$ is simple and injective, $\Cal U(S) = \mod B$ and $\Omega S = M.$ 
Proposition 3.6 asserts that $S$ is semi-Gorenstein-projective.
Also, $S$ is not projective, since $M \neq 0.$ 
Since $S$ is injective and not projective, we see that $S^* = 0,$ thus $S$ is a Nunke module.
$\s$

	\bigskip\bigskip
{\bf 4. (Short) local algebras.}
	\medskip
One may wonder whether there do exist non-trivial 
examples of modules which satisfy the conditions discussed in section 3.
Here we want to mention at least one class of algebras, the short local algebras, were
no examples of this kind do exist. First, let $A$ be an arbitrary local artin algebra.
	\medskip
{\bf 4.1. Proposition.} {\it Let $A$ be a local artin algebra. 
Let $M$ be a semi-Gorenstein-projective module. If $M^*$ is projective, then $M$ is projective.}
	\medskip
Proof. We can assume that $M$ is indecomposable and not projective, thus reduced. 
According to the Main Theorem, there is a minimal complex $P_\bullet = (P_i,f_i)$ 
of projective modules such that $P_\bullet^*$ is acyclic, as exhibited in the display 2.5.
In particular, $M^*$ is the image of $f_0^*$. Since $f_0$ maps into the radical of $P_{-1}$,
$f_0^*$ maps ito the radical of $P_0^*$, thus $M^*$ is a submodule of the radical of a projective 
right module. Since $A$ is local, $M^*$ cannot have an indecomposable projective direct summand.

If we assume that $M^*$ is projective, then $M^* = 0$. But for a local algebra $A$, this
implies that $M = 0,$ a contradiction. This completes the proof.
$\s$
	\bigskip
{\bf 4.2.} We recall that a local algebra $A$ is said to be
{\it short} provided that $J^3 = 0.$ From now on, let $A$ be a short local artin algebra
with $e = |J/J^2|$ and $a = |J^2|$. The pair $(e,a)$ is called the {\it Hilbert type} of $A$.
For any module $M$, we denote by $|M|$ its length and set $t(M) = |M/JM|$. 
If $M$ has Loewy length at most 2, then $\bdim M = (t(M),|JM|)$ is called the {\it dimension
vector} of $M$. Note that if 
$Q_\bullet = (Q_i,d_i\:Q_i\to Q_{i-1})$ is a minimal complex of projective
modules and  $N_i$ be the image of $d_i$, then $N_i$ has Loewy length at most 2, thus its
dimension vector is defined. 
	\bigskip
{\bf  Proposition.} {\it Let $A$ be a short local artin algebra.}
	\smallskip
(a) {\it  Assume that $A$ is not self-injective, and let
$Q_\bullet = (Q_i,d_i\:Q_i\to Q_{i-1})$ be an acyclic minimal complex of projective
modules. Let $N_i$ be the image of $d_i$ and assume that at least one of the modules $N_i$ is
semi-Gorenstein-projective. Then all modules $Q_i$ have the same rank, say rank $t$ and
$\bdim N_j = (t,at)$ for all $j\in \Bbb Z$.}
	\smallskip
(b) {\it 
Let $M$ be a module such that both $M$ and $M^*$ are semi-Gorenstein-projective. 
Then $|\Ker(\phi_M)| = |\Cok(\phi_M)|$.}
(Thus, if $\phi_M$ is a monomorphism or an epimorphism, then $\phi_M$
is an isomorphism, and therefore $M$ is Gorenstein-projective.)
	\medskip
The proof will rely on on two results (Theorems 3 and 4) from 
[RZ3]. 
	\medskip
Proof of (a). Let 
$Q_\bullet = (Q_i,d_i\:Q_i\to Q_{i-1})$ be an acyclic minimal complex of projective
modules, with $N_i$ the image of $d_i$, for all $i\in \Bbb Z$. 
Since $Q_\bullet$ is acyclic and minimal, the canonical maps $Q_i\to N_i$
are projective covers, thus $t(Q_i) = t(N_i)$ for all $i\in \Bbb Z.$  

According to Theorem 3 of [RZ3], the complex $Q_\bullet$ shows that 
$A$ is of Hilbert type $(a+1,a)$ with $a\ge 1$,
and that either all the modules $N_i$ have the same dimension vector (type I), 
in particular all the projective modules $Q_i$ have the same rank, 
or else the rank of the modules $Q_{i}$ strictly increases for $i \gg 0$ (type II). 

Let us now assume that $N_0$ is semi-Gorenstein-projective. Of course, $N_0$ is torsionless
and not projective and
$$
  0 \leftarrow N_0 \leftarrow Q_0 \leftarrow Q_1 \leftarrow \cdots
$$
is a minimal projective resolution of $N_0$. We apply Theorem 4 of [RZ3] to the indecomposable direct summands of $N_0$ and see that $a\ge 2$, 
and that all the modules $N_i = \Omega^i N_0$ with $i\in \Bbb N$ have the same dimension vector
$\bdim N_i = (t,at)$, where $t = t(N_0) = t(Q_0)$. As a consequence, 
$t(Q_i) = t(N_i) = t$ for all $i \ge 0.$ 
This shows that $Q_\bullet$ cannot be of type II. Thus, $Q_\bullet$ is
of type I, and therefore all the projective modules $Q_i$ have the same rank $t$, and  
$\bdim N_i = (t,at)$ for all $i\in \Bbb Z.$ This completes the proof of (a). $\s$
	\medskip 
Proof of (b). If $A$ is self-injective, then all modules are
reflexive, thus (b) holds trivially in this case. We therefore may assume that $A$ 
is not self-injective. 

Let $M$ be a module such that both $M$ and $M^*$ are 
semi-Gorenstein-projective. According to the Main Theorem, there is a minimal complex 
$P_\bullet = (P_i,f_i)$ 
of projective modules such that $P_\bullet^*$ is acyclic, as exhibited in the display 2.5.
We apply 9.3 to the opposite algebra $A^{\op}$, thus to right $A$-modules, namely to the
acyclic complex $P_\bullet^*$ of projective right $A$-modules. Since the image of $f_0^*$
is the semi-Gorenstein-projective module $M^*$, we see that all the modules $P_i^*$ have the
same rank, say $t$. Thus also the modules $P_i$ have rank $t$.

Now $\Omega M$ is the image of $f_1$. Since $P_1$ is a projective cover of $\Omega M$,
we see that $\top \Omega M$ has length $t$. Similarly, $\Tr M^*$ is the image
of $f_{-2}$ and $P_{-2}$ is a projective cover of $\Tr M^*$, thus $\top \Tr M^*$
has length $t$. Next, $\Omega M$ is an idecomposable torsionless semi-projective module 
and not projective, thus Theorem 4 of [RZ3] asserts that its dimension vector is
$(t,at)$. Similarly, $\Tr M^*$ is an indecomposable $\infty$-torsionfree module and not
projective, thus the same reference shows that the dimension vector of $\Tr M^*$ 
is also $(t,at)$. 
We consider the sequence
$$
 0 \to \Omega M @>u>> P_0 @>f_0>> P_{-1} @>f_{-1}>> P_{-2} @>r>> \Tr M^* \to 0,
$$
where $u$ is the canonical inclusion and $r$ the canonical projection. This is a complex,
and the alternating sum of the length of the modules involved is zero 
(the modules $\Omega M$ and $\Tr M^*$ have length $et$, whereas the modules $P_0, P_{-1},P_{-2}$
have length $2et$). It follows that also the alternating sum of the length of the homology
modules is 0, but this is $|\Ker(\phi_M)|-|\Cok(\phi_M)|$.
$\s$
	\bigskip\bigskip
{\bf References.}
	\medskip
\item{[AB]}  M. Auslander and M. Bridger. Stable module theory. Memoirs of the American
  Mathematical Society, No. 94. American Mathematical Society, Providence, R.I., 1969.
\item{[AR]}  M. Auslander, I. Reiten. On a generalized version of the Nakayama conjecture.
   Proc. Amer. Math. Soc. 52 (1975), 69--74.
\item{[ARS]}  M. Auslander, I. Reiten, S. O. Smal\o. Representation Theory
  of Artin Algebras. Cambridge Studies in Advanced Math. 36.
  Cambridge University Press, 1995.
\item{[AM]} L. L. Avramov, A. Martsinkovsky, Absolute, relative, and Tate cohomology of 
  modules of finite Gorenstein dimension, Proc. London Math. Soc. 85(3)(2002), 393–440.
\item{[B]} R.-O. Buchweitz. Maximal Cohen-Macaulay Modules and Tate-Cohomology
  Over Gorenstein Rings. Manuscript, available at http://hdl.handle.net/1807/16682. 1986.
\item{[G]} V. G\'elinas. The depth, the delooping level and the finitistic dimension. 
  Preprint. arXuv:2004.04838v1 
\item{[H]} D. Happel. Homological conjectures in representation theory of finite
   dimensional algebras.
   https://www.math.uni-bielefeld.de/$\sim$sek/dim2/happel2.pdf
\item{[RZ1]}  C. M. Ringel, P. Zhang. Gorenstein-projective and
  semi-Gorenstein-projective modules. 
  Algebra \& Number Theory 14-1 (2020), 1--36.  DOI:
  10.2140/ant.2020.14.1 
\item{[RZ2]}  C. M. Ringel, P. Zhang. Gorenstein-projective and
  semi-Gorenstein-projective modules\. II.
  J. Pure Appl. Algebra 224 (2020) 106248. \newline
  https://doi.org/10.1016/j.jpaa.2019.106248
\item{[RZ3]}  C. M. Ringel, P. Zhang.  Gorenstein-projective modules over short local algebras.
  Preprint.  arXiv:1912.02081 

	\bigskip\bigskip

{\baselineskip=1pt
\rmk
C. M. Ringel\par
Fakult\"at f\"ur Mathematik, Universit\"at Bielefeld \par
POBox 100131, D-33501 Bielefeld, Germany  \par
ringel\@math.uni-bielefeld.de
\smallskip

P. Zhang \par
School of Mathematical Sciences, Shanghai Jiao Tong University \par
Shanghai 200240, P. R. China.\par
pzhang\@sjtu.edu.cn\par}
\bye